\documentclass[12pt]{article}

\makeatletter
\newfont{\footsc}{cmcsc10 at 8truept}
\newfont{\footbf}{cmbx10 at 8truept}
\newfont{\footrm}{cmr10 at 10truept}
\renewcommand{\ps@plain}{%
\renewcommand{\@oddfoot}{\footsc the electronic journal of combinatorics
  {\footbf 10} (2003), \#R00\hfil\footrm\thepage}}
\makeatother
\usepackage{amsmath,amssymb,amsfonts,latexsym,theorem}

\numberwithin{equation}{section}
\newtheorem{theorem}{Theorem}[section]
\newtheorem{proposition}[theorem]{Proposition}
\newtheorem{corollary}[theorem]{Corollary}
\newtheorem{conjecture}[theorem]{Conjecture}
\newtheorem{lemma}[theorem]{Lemma}

{\theorembodyfont\rmfamily
\newtheorem{definition}[theorem]{Definition}
\newtheorem{remark}[theorem]{Remark}
\newtheorem{example}[theorem]{Example}
\newtheorem{notation}[theorem]{Notation}
}

\def\mn{\text{-}}
\def\g{\text{g}}

\def\qed{\hfill $\Box$}

\newenvironment{proof}{\noindent{\scshape Proof.}}{\qed}

\begin{document}

%\pagenumbering{arabic} \pagestyle{headings}
\def\sof{\hfill\rule{2mm}{2mm}}
\def\SS{\mathcal S}
\def\qq{{\bold q}}
\def\mn{\text{-}}

\title{\textbf{Packing patterns into words}}

\author{
\textbf{Alexander Burstein}\\
Department of Mathematics\\
Iowa State University\\
Ames, IA 50011-2064 USA\\
\texttt{burstein@math.iastate.edu} \and
\textbf{Peter H\"ast\"o}\\
Department of Mathematics\\
P.O. Box 4 (Yliopistokatu 5), \\
00014 University of Helsinki, Finland\\
\texttt{peter.hasto@helsinki.fi} \and
\textbf{Toufik Mansour}\\
Department of Mathematics\\
Chalmers University of Technology\\
S-412 96 G\"oteborg, Sweden\\
\texttt{toufik@math.chalmers.se} \and
}

%\date{\today}

\date{\small Submitted: January 1, 2001;  Accepted: January 2, 2001.\\
\small MR Subject Classifications: Primary 05A15; Secondary 05A16}

\maketitle

\begin{abstract}
In this article we generalize packing density problems from
permutations to patterns with repeated letters and generalized
patterns. We are able to find the packing density for some classes
of patterns and several other short patterns.
\end{abstract}

A string $213322$ contains three subsequences $233, 133, 122$ each
of which is \emph{order-isomorphic} (or simply \emph{isomorphic})
to the string $122$, i.e. ordered in the same way as $122$. In
this situation we call the string $122$ a \emph{pattern}.

Herb Wilf first proposed the systematic study of pattern
containment in his 1992 address to the SIAM meeting on Discrete
Mathematics. However, several earlier results on pattern
containment exist, for example, those by Knuth \cite{Knuth} and
Tarjan \cite{Tarjan}.

Most results on pattern containment actually deal with
\emph{pattern avoidance}, in other words, enumerate or consider
properties of strings over a totally ordered alphabet which avoid
a given pattern or set of patterns.

There is considerably less research on other aspects of pattern
containment, specifically, on packing patterns into strings over a
totally ordered alphabet (but see \cite{AAHHS, BSV, Hasto, Price,
Stromquist}). In fact, all pattern packing except the one in
\cite{Stromquist} (later generalized in \cite{AAHHS}) dealt with
packing permutation patterns into permutations (i.e. strings
without repeated letters). In this paper, we generalize the
packing statistics and results to patterns over strings with
repeated letters and relate them to the corresponding results on
permutations.

\section{Preliminaries} \label{sec:prelim}

Let $[k]=\{1,2,\dots,k\}$ be our canonical totally ordered
alphabet on $k$ letters, and consider the set $[k]^n$ of
$n$-letter words over $[k]$. We say that a pattern $\pi\in[l]^m$
\emph{occurs} in $\sigma\in[k]^n$, or $\pi$ \emph{hits} $\sigma$,
or that $\sigma$ \emph{contains} the pattern $\pi$, if there is a
subsequence of $\sigma$ order-isomorphic to $\pi$.

Given a word $\sigma\in[k]^n$ and a set of patterns
$\Pi\subseteq[l]^m$, let $\nu(\Pi,\sigma)$ be the total number of
occurrences of patterns in $\Pi$ ($\Pi$-patterns, for short) in
$\sigma$. Obviously, the largest possible number of
$\Pi$-occurrences in $\sigma$ is $\binom{n}{m}$, when each
subsequence of length $m$ of $\sigma$ is an occurrence of a
$\Pi$-pattern. Define
\[
\begin{split}
\mu(\Pi,k,n)&=\max\{\,\nu(\Pi,\sigma)\,|\,\sigma\in[k]^n\},\\
d(\Pi,\sigma)&=\frac{\nu(\Pi,\sigma)}{\binom{n}{m}}\ {\rm and}\\
\delta(\Pi,k,n)&=\frac{\mu(\Pi,k,n)}{\binom{n}{m}}=
\max\{\,d(\Pi,\sigma)\,|\,\sigma\in[k]^n\},
\end{split}
\]
respectively, the maximum number of $\Pi$-patterns in a word in
$[k]^n$, the probability that a subsequence of $\sigma$ of length
$m$ is an occurrence of a $\Pi$-pattern, and the maximum such
probability over words in $[k]^n$. We want to consider the
asymptotic behavior of $\delta(\Pi,k,n)$ as $n\to\infty$ and
$k\to\infty$.

\begin{proposition} \label{prop:monotone}
If $n>m$, then $\delta(\Pi,k,n)\le\delta(\Pi,k,n-1)$ and
$\delta(\Pi,k,n)\ge\delta(\Pi,k-1,n)$.
\end{proposition}

\begin{proof}
The proof of Proposition~1.1 in \cite{AAHHS} also applies to the
first inequality in our proposition, since possible
repetition of letters is irrelevant here. To see that the second
inequality is true, note that increasing $k$, i.e. allowing more
letters in our alphabet, can only increase $\mu(\Pi,k,n)$, and
hence $\delta(\Pi,k,n)$.
\end{proof}

The greatest possible number of distinct letters in a word
$\sigma$ of length $n$ is $n$, which implies that
$\mu(\Pi,k,n)=\mu(\Pi,n,n)$ for $k\ge n$, and hence,
$\delta(\Pi,k,n)=\delta(\Pi,n,n)$ for $k\ge n$. Therefore,
\[
\delta(\Pi,n,n)=\lim_{k\to\infty}{\delta(\Pi,k,n)}.
\]
We also have
$\delta(\Pi,n,n)=\delta(\Pi,n+1,n)\ge\delta(\Pi,n+1,n+1)$, so
$\delta(\Pi,n,n)$ is non-increasing and nonnegative, and there
exists
\[
\delta(\Pi)=\lim_{n\to\infty}{\delta(\Pi,n,n)}
=\lim_{n\to\infty}{\lim_{k\to\infty}{\delta(\Pi,k,n)}}.
\]
We call $\delta(\Pi)$ the \emph{packing density} of $\Pi$.

Obviously, there are two double limits. Since
$0\le\delta(\Pi,k,n)\le 1$, it immediately follows that there
exists
\[
\delta(\Pi,k)=\lim_{n\to\infty}{\delta(\Pi,k,n)}\in[0,1]
\]
and that $\{\delta(\Pi,k)\,|\,k\in\mathbb{N}\}$ is nondecreasing
as $k\to\infty$. Hence, there exists
\[
\delta'(\Pi)=\lim_{k\to\infty}{\delta(\Pi,k)}
=\lim_{k\to\infty}{\lim_{n\to\infty}{\delta(\Pi,k,n)}}.
\]

It is easy to see that $\delta'(\Pi)\le\delta(\Pi)$. Naturally,
one wishes to determine when $\delta'(\Pi)=\delta(\Pi)$. In this
paper, we will provide a sufficient condition for this equality.

The set $[k]^n$ is finite, so for each $k$ and $n$, there is a
string $\sigma(\Pi,k,n)\in[k]^n$ such that
$d(\Pi,\sigma(\Pi,k,n))=\delta(\Pi,k,n)$. To find $\delta(\Pi)$,
we will need to find $\delta(\Pi,k,n)$, hence maximal
$\Pi$-containing permutations $\sigma(\Pi,k,n)$ are of interest to
us, especially, their asymptotic shape as $n\to\infty$ and
$k\to\infty$.

\begin{example} \label{ex:1...1}
Let $\Pi=\{c_m\}$, where $c_m$ is a constant string of $m$ 1's.
Then, clearly, $\sigma(\Pi,k,n)=c_n$ and $d(c_m,c_n)=1$ for $n\ge
m$, so $\delta(c_m,k,n)=1$ for $n\ge m$, and hence
$\delta'(c_m)=\delta(c_m)=1$ for any $m\ge 1$.
\end{example}

\begin{example} \label{ex:id_m}
Let $\Pi=\{id_m\}$, where $id_m$ is the identity permutation of
$S_m$. Then $\sigma(id_m,n,n)=id_n$, so $d(id_m,id_n)=1$,
$\delta(id_m,n,n)=1$ and $\delta(id_m)=1$.

Determining $\delta'(id_m)$ is a bit harder. It is easy to see
that $\sigma(id_m,k,n)$ must be a nondecreasing string of digits
in $[k]$. Let $n_i$ be the number of digits $i$ in
$\sigma(id_m,k,n)$, then
$\mu(id_m,k,n)=\nu(id_m,\sigma(id_m,k,n))=n_1n_2\dots n_k$ and
$n_1+n_2+\dots+n_k=n$. To maximize the above product we need
$n_1=n_2=\dots=n_k=\frac{n}{k}$. (More exactly, \cite{Price} shows
that we should choose for $n_i$'s to be such integers that
$|n_i-\frac{n}{k}|<1$ and $|n_1+\dots+n_r-\frac{rn}{k}|<1$ for
each $r=1,2,\dots,k$.) It follows that
\[
\delta(id_m,k,n)\sim
\frac{\binom{k}{m}\left(\frac{n}{k}\right)^m}{\binom{n}{m}}
\]
(where $a_n\sim b_n$ means $\lim_{n\to\infty}{a_n/b_n}=1$), so
$\delta(id_m,k)= \binom{k}{m}\frac{m!}{k^m}$, and thus
$\delta'(id_m)=1$ as expected.
\end{example}

Packing density was initially defined for patterns in
permutations. Therefore, we must show that the packing density on
permutations agrees with the packing density on words.

\begin{theorem} \label{thm:perm}
Let $\Pi\subseteq S_m$ be a set of permutation patterns, then
\[
\delta(\Pi)=\lim_{n\to\infty}{\frac{\max\{\,\nu(\Pi,\sigma)\,|\,\sigma\in
S_n\}}{\binom{n}{m}}},
\]
i.e. the packing density of $\Pi$ on words is equal to that on
permutations.
\end{theorem}
\begin{proof}
It is enough to prove that
\[
\mu(\Pi,n,n)=\max\{\,\nu(\Pi,\sigma)\,|\,\sigma\in S_n\},
\]
in other words, that there is a permutation in $S_n$ among the
maximal $\Pi$-containing words in $[n]^n$. Consider any maximal
$\Pi$-containing word $\sigma\in[n]^n$. Let $n_i$ be the
multiplicity of the letter $i$ in $\sigma$. Let $i_j$ denote the
$j$th occurrence of the letter $i$, and consider the map
$f:[n]^n\to S_n$ induced by the map
$i_j\mapsto\sum_{r=1}^{i}{n_r}-j+1$. Since all letters of each
pattern in $\Pi$ are distinct, $\Pi$ occurs in $f(\sigma)$ at
least at the same positions $\Pi$ occurs in $\sigma$, so
$\nu(\Pi,f(\sigma))\ge\nu(\Pi,\sigma)$. The rest is easy.
\end{proof}

Apart from computing packing densities of patterns, we would also
like to determine which patterns have equal packing densities,
which ones are asymptotically more packable than others, etc. For
example, it is easy to see that the packing density is invariant
under the usual symmetry operations on $[l]^m$: \emph{reversal}
$r:\tau(i)\to \tau(m-i+1)$ and \emph{complement} $c:\tau(i)\to
l-\tau(i)+1$, (packing density is also invariant under inverse
$i:\tau\to\tau^{-1}$ when packing permutations into permutations).
The operations $r$ and $c$ generate $D_2$, while $r,c,i$ generate
$D_4$. Patterns which can be obtained from each other by a
sequence of symmetry operations are said to belong to the same
\emph{symmetry class}.

\begin{example} \label{ex:pd-S3}
The symmetry class representatives of patterns in $[3]^3$ are
$111$, $112$, $121$, $123$ and $132$. We know that
$\delta(111)=1=\delta(123)$. Galvin, Kleitmann and Stromquist
(independently, unpublished, see chronology in \cite{Price})
showed that $\delta(132)=2\sqrt{3}-3\approx 0.4641$. Thus, we only
need to determine the packing densities of 112 and 121 to
completely classify patterns of length 3.
\end{example}

Price \cite{Price} extended Stromquist's results \cite{Stromquist}
to packing a single pattern $\pi=1m(m-1)\dots2$ and handled other
single patterns such as $2143$. Since we will also be concerned
mostly with singleton sets of patterns $\Pi=\{\pi\}$, we will
write $\delta(\pi)$ for $\delta(\{\pi\})$, etc.

Price's results deal with patterns of specific type, the so-called
\emph{layered} patterns.

\begin{definition} \label{def:layered}
A \emph{layered} pattern is a strictly increasing sequence of
strictly decreasing substrings. These substrings are called the
\emph{layers} of $\sigma$.
\end{definition}

\begin{notation}
It easy to see that a layered pattern is uniquely determined by
the sequence of its layer lengths, hence we may denote it by such
sequence, e.g.
$\widehat{321}\,\widehat{54}\,\widehat{9876}=[3,2,4]$,
$\widehat{1}\widehat{2}\widehat{3}=[1,1,1],
\widehat{1}\widehat{32}=[1,2],\widehat{21}\widehat{3}=[2,1],\widehat{321}=[3]$
are layered, with layers denoted by hats, while $312, 231$ are
non-layered.
\end{notation}

In fact, note that the union of symmetry classes of layered
patterns consists of exactly the permutations avoiding patterns in
the symmetry classes of $1342, 1423, 2413$.

In \cite{Stromquist}, Stromquist proved a theorem (later
generalized in \cite{AAHHS}) on packing layered patterns into
permutations.
The inductive proof of this theorem defines a permutation (or a
poset) $\pi$ to be \emph{layered on top} (or \emph{LOT}) if any of
its maximal elements is greater than any non-maximal element. The
set of these maximal elements is called the \emph{final layer} of
$\pi$ (even if $\pi$ is not necessarily layered).
\begin{proposition} \label{prop:lot}
Let $\Pi$ be a multiset of LOT permutations (not necessarily all
distinct or of equal length). Then there is an LOT permutation
$\sigma^{*}$ which maximizes the expression
\begin{equation} \label{eq:lot}
\nu(\Pi,\sigma)=\sum_{\pi\in\Pi}{a_{\pi}\nu(\pi,\sigma)}, \quad
a_\pi\ge 0.
\end{equation}
Furthermore, if the final layer of every $\pi\in\Pi$ has size
greater than 1, then every such $\sigma^{*}$ is LOT.
\end{proposition}

Applying this proposition inductively, \cite{AAHHS}, following
\cite{Stromquist}, obtains
\begin{theorem} \label{thm:layered}
Let $\Pi$ be a multiset of layered permutations. Then there is a
layered permutation $\sigma^{*}$ which maximizes the expression
(\ref{eq:lot}). Furthermore, if all the layers of every
$\pi\in\Pi$ have size greater than 1, then every such $\sigma^{*}$
is layered.
\end{theorem}

Following \cite{AAHHS,Price}, we will also define the
\emph{$\ell$-layer packing density} $\delta_\ell(\Pi)$ for sets of
layered permutations $\Pi$ as the packing density of $\Pi$ among
the permutations with at most $\ell$ layers. It was shown in both
of the above works that
$\displaystyle\delta(\Pi)=\lim_{\ell\to\infty}{\delta_\ell(\Pi)}$.

\section{Monotone and layered patterns} \label{sec:monotone}

The easiest type of patterns with repeated letters are those whose
letters are nondecreasing (or non-increasing) from left to right.
By analogy with layered patterns, we will consider nondecreasing
patterns.

We will call a maximal constant segment of a word a \emph{block}.
For a letter $a$ and integer $k\ge 1$, we will define
$a^k=\underbrace{a\dots a}_k$.

\begin{theorem} \label{thm:monotone}
Let $\Pi\in[l]^m$ be a set of nondecreasing patterns
$\pi=1^{a_1(\pi)}2^{a_2(\pi)}\dots l^{a_l(\pi)}$. For each
$\pi\in\Pi\subseteq[l]^m$, let $\hat\pi\in S_m$ be the layered
pattern $\hat\pi=[a_1(\pi),\dots,a_l(\pi)]$, and let
$\hat\Pi=\{\hat\pi\,|\,\pi\in\Pi\}$. Then
$\delta(\Pi,k)=\delta_k(\hat\Pi)$ and
$\delta'(\Pi)=\delta(\Pi)=\delta(\hat\Pi)$.
\end{theorem}

\begin{proof}
There is a natural bijection between nondecreasing patterns on $l$
letters and layered patterns with $l$ layers. The map $f$ of
Theorem~\ref{thm:perm}, induced by the map
$i_j\mapsto\sum_{r=1}^{i}{a_r(\pi)}-j+1$ (where $i_j$ is the $j$th
$i$ from the left), maps $\pi$ to $\hat\pi$. Clearly, $f^{-1}$ is
induced by a map which takes each element in the $i$th layer (the
$i$th basic subsequence, in general) to integer $i$.
\end{proof}

\begin{example} \label{ex:pd-112-1122}
Using the previous theorem and results of Price \cite{Price}, we obtain
$\delta(112)=\delta(\widehat{21}\widehat{3})=2\sqrt{3}-3$,
$\delta(1122)=\delta(\widehat{21}\widehat{43})=3/8$. More
generally, for $k\ge2$,
\[
\delta(\underbrace{1\dots1}_k 2)=ka(1-a)^{k-1},\ \mbox{ where } \
0<a<1,\ \ ka^{k+1}-(k+1)a+1=0.
\]
Similarly, for $r,s\ge 2$,
\[
\delta(\underbrace{1\dots1}_r\underbrace{2\dots2}_s)
=\delta(\underbrace{1\dots1}_r\underbrace{2\dots2}_s,\,2)
=\binom{r+s}{r}\frac{r^rs^s}{(r+s)^{r+s}}.
\]
Using the results of Albert et al. \cite{AAHHS}, we also find that
$\delta(1123)=\delta(1233)=\delta(1243)=3/8$,
$\delta(\{122,112\})=\delta(\{132,213\})=3/4$.
\end{example}

\begin{notation}
A monotone nondecreasing pattern is uniquely determined by the
sequence of its block lengths. Because of this and as a
consequence of Theorem~\ref{thm:monotone}, we may by abuse of
notation denote a monotone nondecreasing pattern by the sequence
of its block lengths, e.g. $112=[2,1],122=[1,2],123=[1,1,1]$.
\end{notation}

By analogy with layered permutations, we define layered strings as
follows.

\begin{definition} \label{def:weakly}
A string $\pi\in[l]^m$ is \emph{layered} if it is a
concatenation of a strictly increasing sequence of non-increasing
substrings. In other words, $\pi=\pi_1\dots\pi_r$, where $\pi_i$
are non-increasing, and $\pi_1<\dots<\pi_r$ (that is any
letter of $\pi_i$ is less than any letter of $\pi_j$ if
$i\le j$). Substrings $\pi_i$ maximal with respect to these
properties are called the layers of $\pi$.
\end{definition}

\begin{definition}
Let us say that the layered permutation $\pi$ is
\textit{simple} if there exists a sequence $\{\sigma_n\}$ of layered
permutations with $\sigma_n\in S_n$
such that every $\sigma_n$ has $r$ layers and $\lim_{n\to\infty} d(\pi,\sigma_n)
= \delta(\pi)$.
\end{definition}

Simple permutations are, as indicated by the name, the easiest type of
permutations to calculate the packing density of. Indeed, it was show in
\cite[Theorem~1.2]{Hasto} that the layered permutation $\pi$ of type
$[m_1, \ldots, m_r]$ with $\log_2(r+1) \le \min\{m_i\}$ is simple
and that in this case
$$ d(\pi) = \frac{m!}{m^m} \prod_{k=1}^r \frac{m_k^{m_k}}{m_k!}, $$
where $m:=m_1+\ldots +m_r$.

\begin{theorem}\label{thm:reduction}
Let $\pi$ be layered pattern with each layer isomorphic to either
$k\ldots 1$ or $1\ldots 1$. Let $\pi'$ be the layered permutation
with layer lengths equal to those of $\pi$. If $\pi$ is simple,
then $\delta(\pi)=\delta(\pi')$.
\end{theorem}

\begin{proof}
Let us denote by $m$ the number of layers in $\pi$.
 If $f$ is an operation as in Theorems \ref{thm:perm} and
\ref{thm:monotone} and $\pi$ is layered, then $\pi'=f(\pi)$ is
layered. Since $\pi'$
 is simple, the $\pi'$-maximal permutation is essentially one
 with $m$ layers of size proportional to those of $\pi$. But
 transforming this permutation into a layered pattern by changing
 a layer to block if the corresponding layer of $\pi$ is a block
 gives a pattern $\sigma$ for which $d(\pi, \sigma)\to \delta(\pi')$.
Therefore $\delta(\pi)\ge \delta(\pi')$.

Let $\sigma$ be a $\pi$-maximal pattern. Then every occurence
of $\pi$ in $\sigma$ is an occurence of $f(\pi) = \pi'$ in $f(\sigma)$,
so that $\delta(\pi) \le \delta(\pi')$. It follows that
$\delta(\pi) = \delta(\pi')$.
\end{proof}

\begin{example}\label{thm:k+1}
Let $\pi = k(k-1) \ldots 1 (k+1)^q$. If $q\ge 2$ then
\[
\delta(\pi) = \binom{k+q}{k} \frac{k^k q^q}{(k+q)^{k+q}}.
\]
If $q=1$ then $ \delta(\pi) = \delta(1^k 2)=\delta([k,1])$, given
in Example \ref{ex:pd-112-1122}.
The first claim follows by Theorem~\ref{thm:reduction} and
\cite[Theorem~1.2]{Hasto}. The second claim follows by
Theorem~\ref{thm:reduction} and \cite[Theorem~5.2]{Price}.
\end{example}

\begin{conjecture} \label{conj:weakly}
If $\Pi$ is a set of %weakly
layered patterns, then
$\delta'(\Pi)=\delta(\Pi)$ and among maximal $\Pi$-containing
strings in $[k]^n$, there is one which is %weakly
layered.
\end{conjecture}

Next we will discuss a non-monotone type of patterns related to
monotone patterns.

\begin{theorem}\label{thm:pqr} Let $\pi = 1^p 2^r 1^q$, for $p,q,r\ge 1$. Then
\[
\delta(\pi)=\binom{p+q}{p}\frac{p^p q^q}{(p+q)^{p+q}}\,
\delta(1^{p+q}2^r)=\binom{p+q}{p}\frac{p^p q^q}{(p+q)^{p+q}}\,
\delta([p+q,r]).
\]
\end{theorem}

\begin{proof} Let $\sigma$ be a $\pi$-maximal pattern of length $n$.
Denote by $a_i$ the number of $i$'s in $\sigma$. It is clear that
$\sigma$ can be assumed to have at least two blocks at every
height except the greatest.

Let us compare the hits (occurrences) of $\pi$ in $\sigma$ with
those of the pattern $\pi' = 1^{p+q} 2^r$ in $\sigma' = 1^{a_1}
2^{a_2} \ldots k^{a_k}$. Lets count the number of hits in each
case with the blocks of $1$'s at height $i$ and the $2$'s at
height $j>i$. The maximum number of such hits of $\pi$ in $\sigma$ occurs in
the pattern $i^{pa_i/(p+q)} j^{a_j} i^{qa_i/(p+q)}$ and equals
$$ \binom{pa_i/(p+q)}{p}  \binom{a_j} {r} \binom{qa_i/(p+q)}{q}. $$
(This argument is strictly true only if $a_i$ is divisible by
$p+q$, otherwise we have to round suitably.) On the other hand the
hits of $\pi'$ in $\sigma'$ with the $1$'s and the $2$'s at these
heights occurs in
\[
\binom{a_i}{p+q} \binom{a_j}{r}
\]
cases. By considering this ratio for large $a_i$, we find that
\[
\nu (\pi, \sigma) \le \binom{p+q}{p} \frac{p^p q^q}{(p+q)^{p+q}}
\nu (\pi', \sigma').
\]
(We do not need to consider small $a_i$'s since their contribution
as $n\to \infty$ will be negligible.) But we know the density of
$\pi'$ by Theorem~\ref{thm:monotone}, and so it follows that
\[
\delta(\pi) \le \binom{p+q}{p} \frac{p^p q^q}{(p+q)^{p+q}}
\delta([p+q, r]).
\]

On the other hand it is easy to see that we can construct patterns
containing this many $\pi$'s; we take a $\pi'$-maximal pattern and
split each block except the one on the highest level into two
blocks of relative sizes $p$ and $q$ and place the first before
and the latter after all higher height blocks. Therefore the
inequality is in fact is an equality, and the theorem is proved.
\end{proof}

\begin{remark} \label{rem:r2}
If $r>1$ in the previous theorem, then
\[
\delta([p+q, r])= {p+q+r \choose r} \frac{(p+q)^{p+q} r^r}
{(p+q+r)^{p+q+r}},
\]
and so
\[
\delta(\pi) = {p+q+r \choose p,q,r} \frac{p^p q^q r^r}
{(p+q+r)^{p+q+r}}.
\]
The $\pi$-maximizing string here is of the type $1^a 2^c 1^b$ with
asymptotic layer lengths
\[
\left(\frac{p}{p+q+r},\frac{r}{p+q+r},\frac{q}{p+q+r}\right).
\]
\end{remark}

\begin{remark} \label{r1}
When $r=1$, we can calculate $\delta([p+q, r])$ as in Example
\ref{ex:pd-112-1122}, which yields
\[
\delta(\pi)=\binom{p+q}{p} p^p q^q
\left(1-(p+q)\alpha\right)\,\alpha^{p+q-1},
\]
where $\alpha\in(0,1)$ is the unique solution of
$(1-sx)^{s+1}=1-(s+1)x$ and $s=p+q$. It is easy to see that
\[
\alpha=\frac{1}{s+1}-(s+1)^{-(s+2)}+O\left((s+1)^{-2s}\right),
\]
since for $x_0 = 1/(s+1) - (s+1)^{-(s+2)}$ we have
\[
\begin{split}
(1-sx_0)^{s+1} + (s+1)x_0 - 1 &= \left(\frac{1}{s+1} +
\frac{s}{(s+1)^{s+2}}\right)^{s+1} - (s+1)^{-(s+1)} \\
&= \frac{s}{(s+1)^{2s+1}} + O\left((s+1)^{1-3s}\right)
\end{split}
\]
For $s\ge 3$, the error in $\alpha$ is at most $4^{-6}<0.00025$,
so $x_0$ approximates $\alpha$ up to at least 3 decimal places.
Note that Theorem~\ref{thm:reduction} also applies when $p=0$ or
$q=0$. Note also that for $a$ in Example \ref{ex:pd-112-1122} we
have $a=1-s\alpha$.

The $\pi$-maximizing string here is of the type
$1^{a_1}2^{a_2}3^{a_3}\dots\dots3^{b_3}2^{b_2}1^{b_1}$ with
asymptotic layer lengths
$(p\alpha,p\alpha(1-s\alpha),p\alpha(1-s\alpha)^2,\dots,
\dots,q\alpha(1-s\alpha)^2,q\alpha(1-s\alpha),q\alpha)$.
\end{remark}

\begin{example} \label{thm:121}
$\delta(121)=\frac{1}{2}\delta(112)=\frac{1}{2}\delta(213)=\sqrt{3}-3/2$.
This completes the inventory of
packing densities of 3-letter patterns by symmetry class.
\begin{center}
\begin{tabular}{|c||c|c|c|c|c|} \hline
  % after \\: \hline or \cline{col1-col2} \cline{col3-col4} ...
  Symmetry class & 111 & 112 & 121 & 123 & 132 \\ \hline
  Packing density & 1 & $2\sqrt{3}-3$ & $\displaystyle\frac{2\sqrt{3}-3}{2}$ & 1
& $2\sqrt{3}-3$ \\ \hline
\end{tabular}
\end{center}
\end{example}

\section{Generalized patterns} \label{sec:generalized}

Here we consider packing generalized patterns into words.
\emph{Generalized patterns} were introduced by Babson and
Steingr\'{\i}msson~\cite{BS} and allow the requirement that some
adjacent letters in a pattern be adjacent in its occurrences in an
ambient string as well. For example, an occurrence of a
generalized pattern $21\mn3$ in a permutation $\pi = a_1 a_2
\cdots a_n$ is a subsequence $a_i a_{i+1} a_j$ of $\pi$ such that
$a_{i+1}<a_i<a_j$. Clearly, in the new notation, classical
patterns are those with all hyphens, such as $1\mn3\mn2$.

\begin{notation}
This notation (introduced in \cite{BS}) may be a little confusing
since classical patterns (the ones with all hyphens) were
previously written the same way as the generalized patterns with
all adjacent letters (i.e. with no hyphens). From now on, we will
use the generalized pattern notation. However, if we consider
subword patterns (those with no hyphens), we may write $\pi_\g$
for a generalized pattern $\pi$ without hyphens where the context
allows for ambiguity.
\end{notation}

If $\pi\in[l]^m$ is a generalized pattern with $b$ blocks of
consecutive letters (i.e. $b-1$ hyphens), then it is easy to see
by considering the positions of the first letters of the blocks of
$\pi$ that the maximum possible number of times $\pi$ can occur in
$\sigma\in[k]^n$ is at most
\[
\binom{n-m+b}{b}
\]
(this yields $\binom{n}{m}$ when $b=m$, i.e. when $\pi$ is a
classical pattern).

In fact, this maximum is achieved when $\pi$ is a \emph{constant}
generalized pattern, i.e. any of the generalized patterns obtained
from the constant strings $11\dots1$ by inserting hyphens at
arbitrary positions (possibly, none). Obviously, maximal
$\pi$-containing strings are the constant strings of length $n$.
Thus, any set of constant generalized patterns has packing density
1. Similarly, any set $\Pi$ of hyphenated identity generalized
patterns has $\delta(\Pi)=1$.

Given a set of generalized patterns with $b$ blocks,
$\Pi\subseteq[l]^m$, we define the packing density of $\Pi$
similarly to that of a set of classical patterns. We will use the
same notation as in Section \ref{sec:prelim} for the generalized
patterns.

It is not difficult to see that the analog of Theorem
\ref{thm:perm} holds for generalized patterns as well.
\begin{theorem} \label{thm:gen-perm}
Let $\Pi\subseteq S_m$ be a set of generalized permutation
patterns, then the packing density of $\Pi$ on words is equal to
that on permutations.
\end{theorem}
\begin{proof}
The same argument as in Theorem~\ref{thm:perm} shows that among
maximal $\Pi$-containing strings in $[n]^n$ there is one that has
no repeated letters.
\end{proof}

\subsection{Generalized patterns without hyphens}

\begin{theorem} \label{thm:gen-monotone}
Let $\pi=1^{a_1}2^{a_2}\dots l^{a_l}\in[l]^m \ (l>1)$ be a
nonconstant monotone generalized pattern without hyphens.
If there exists a positive integer $j\le l-2$ such that
$a_1\le a_{j+1}$, $a_{i}=a_{i+j}\ (2\le i\le l-j-1)$ and $a_{l-j}\ge a_l$,
then we denote by $j_0$ the least such $j$ and define
$M_\pi = a_2+\dots+a_{j_0+1}$. Otherwise we set
$M_\pi = \max\{a_1, a_l\} + a_2 + \ldots + a_{l-1}$.
In either case we have $\delta(\pi)=\delta'(\pi)=1/M_\pi$.
\end{theorem}

\begin{proof}
$M_\pi$ is the smallest shift at which $\pi$ overlaps with itself.
The rest is clear.
\end{proof}

\begin{theorem} \label{thm:gen-layered}
Let $\pi=[a_1,a_2,\dots,a_l]\in S_m$ be any $l$-layer $(l>1)$
generalized pattern without hyphens. Let $M_\pi$ be as in Theorem
\ref{thm:gen-monotone}. Then $\delta(\pi)=\delta'(\pi)=1/M_\pi$.
\end{theorem}
\begin{proof}
The same mapping as in Theorem~\ref{thm:monotone} shows that our
$\pi$ has the same packing density as the corresponding monotone
generalized pattern without hyphens of Theorem
\ref{thm:gen-monotone}.
\end{proof}

\begin{corollary}
Let $\pi_1=11\dots12_g\in[2]^m$ and
$\pi_2=1m(m-1)\ldots2_g\in[m]^m$, then
$\delta(\pi_1)=\delta'(\pi_1)=1/(m-1)$ and
$\delta(\pi_2)=\delta'(\pi_2)=1/(m-1)$.
\end{corollary}

For instance, $\delta(112_g)=\delta'(112_g)=1/2$,
$\delta(132_g)=\delta'(132_g)=1/2$ and
$\delta(123_g)=\delta'(123_g)=1$.

\subsection{Generalized patterns with one hyphen}

The maximal number of occurrences of a generalized pattern in
$[l]^m$ with one hyphen (i.e. with $b=2$ blocks) is
$\binom{n-m+2}{2}\sim n^2/2$ as $n\to\infty$.

\begin{proposition} \label{prop:11-2}
$\delta(11\mn2)=\delta'(11\mn2)=1$.
\end{proposition}
\begin{proof}
Let $\sigma\in [k]^n$ be a maximal $(11\mn2)$-containing word,
then $\sigma$ is a monotone nondecreasing string in which letter
$i$ occurs $n_i$ times, $n_1+\cdots+n_k=n$. Then
$\mu(11\mn2,n,k)=\max\{\sum_{i=1}^k {(n_i-1)(n_{i+1}+\cdots+n_k)}\
:\ n_1+\dots+n_k=n\}$. From here, it is not difficult to determine
that $\mu(11\mn2,n,k)\sim n^2/2$ as $n\to\infty$. Choose $n_i$'s
to be such integers that $|n_i-\frac{n}{k}|<1$ and
$|n_1+\dots+n_r-\frac{rn}{k}|<1$ for each $r=1,2,\dots,k$. Then
\[
\mu(11\mn2,n,k)\sim \left(\frac{n}{k}\right)^2\binom{k}{2},
\]
out of $\binom{n-1}{2}$ maximum possible occurrences, and the
result follows.
\end{proof}

%\begin{proposition} \label{prop:12-1}
%$\delta(12\mn1)=\delta'(12\mn1)=1/3$.
%\end{proposition}
%\begin{proof}
%Let $\sigma\in [k]^n$ be a word with maximum occurrences of
%$12\mn1$, then
%$$ \sigma=1212\cdots1211..1\in[2]^n $$
%where the string $12$ occurs in $\sigma$ exactly $d$ times. So
%$$ \mu(12\mn1,n,k)=\max_{1\le d\le n}(d(d-1)/2+d(n-2d)), $$
%and the maximum occurs at $d\sim n/3$. The rest is easy to check.
%\end{proof}

\begin{proposition} \label{prop:perm1dash}
$\delta(12\mn3)=\delta(21\mn3)=1$.
\end{proposition}
\begin{proof}
For pattern $12\mn3$, consider the identity permutation. For
pattern $21\mn3$, consider the layered permutations of length $n$
with $\sqrt{n}$ layers of length $\sqrt{n}$.
\end{proof}

\smallskip

We think, but have not been able to prove rigorously, that
$\delta(12\mn1)=\delta'(12\mn1)=1/3$. At least $\delta(12\mn1, 2)
= 1/3$, since in this case the string with the maximal number of
occurrences of $12\mn1$ is of the type
$$ \sigma=1212\cdots1211..1\in[2]^n $$
where the string $12$ occurs in $\sigma$ exactly $d$ times. So
$$ \mu(12\mn1,n,2)=\max_{1\le d\le n}(d(d-1)/2+d(n-2d)), $$
and the maximum occurs at $d\sim n/3$. It seems that allowing
more symbols in $\sigma$ does not change anything, but here
we could not find a proof.

A more general question related to this and somewhat analogous to
the question of simple layered permutations is: for which $\pi\in
[k]^n$ is $\delta(\pi,k)=\delta(\pi)$?

\section{The problem of the shortest common superpattern}

This problem deals with packing different patterns into a word.
Let $n(l,m)$ be the length of the shortest word which contains
every pattern of length $m$ on at most $l$ letters. Clearly,
$n(l,m)=n(m,m)$ for $m\le l$, hence we are interested only in the
values of $n(l,m)$ for $m\ge l$.

For example, $n(2,2)=3$ (since $121$ contains patterns $11,12,21$)
and $n(3,3)=7$ (since $1231231$ contains patterns
$111,112,121,211,122,212,221,123,132,213,231,312,321$).

\begin{lemma}
For $m\ge l$, $n(l,m)\le l(m-1)+1$.
\end{lemma}
\begin{proof}
Consider the word $\tau=(id_l)^{m-1}1$ where $id_l=123\dots l$.
The rest is obvious.
\end{proof}

\smallskip

At least in the case of $n(l,l)$, this upper bound is apparently a
lower bound as well, although we have not been able to prove it.

\begin{conjecture}
For any $l\ge 1$, $n(l,l)=l^2-l+1$.
\end{conjecture}

This differs from the corresponding result in \cite{EELW} on
permutation patterns, i.e. those in $S_l$, where the upper bound
of $3l^2/4$ for the length of the shortest common superpattern was
established, and there is numerical evidence that the actual value
is closer to $l^2/2$.

%======================References========================

\end{document}